\documentclass[11pt,twoside]{article}
\usepackage{latexsym}
\usepackage{amssymb,amsbsy,amsmath,amsfonts,amssymb,amscd}
\usepackage{graphicx}
\usepackage{hyperref}
\usepackage{xcolor}
\newcommand{\commentout}[1]{}

\newcommand{\R}{\mathbb{R}}

\newcommand {\e}  {\varepsilon}

\newcommand {\lb} {\lambda}
\newcommand {\Chi} {{\bf \raise 2pt \hbox{$\chi$}} }
\newcommand {\f}   {\frac}
\newcommand {\p}   {\partial}
\newcommand{\dis}{\displaystyle}
\newcommand {\proof} {\noindent {\bf Proof}. }
\newcommand{\beq}{\begin{equation}}
\newcommand{\eeq}{\end{equation}}
\newcommand{\bea} {\begin{array}{rl}}
\newcommand{\eea} {\end{array}}
\newcommand{\bepa}{\left\{ \begin{array}{l}}
\newcommand{\eepa} {\end{array}\right.}
\newtheorem{theorem}{Theorem}[section]
\newtheorem{lemma}[theorem]{Lemma}

\newcommand {\no}{\noindent}
\newcommand{\qed}{{ \hfill
                       {\unskip\kern 6pt\penalty 500 \raise -2pt\hbox{\vrule\vbox to 6pt{\hrule width 6pt
                       \vfill\hrule}\vrule} \par}   }}
\title{ {Rare mutations limit of a steady state dispersal evolution model }}

 \author{
Beno\^ \i t Perthame\thanks {Sorbonne Universit\'es, UPMC Univ Paris 06, CNRS, INRIA,  UMR 7598, Laboratoire Jacques-Louis Lions, \'Equipe MAMBA, 4, place Jussieu 75005, Paris, France. Email:~benoit.perthame@upmc.fr} 
\and Panagiotis E. Souganidis\thanks{Department of Mathematics, University of Chicago, Chicago, IL 60637, USA. Email: souganidis@math.uchicago.edu}
\thanks {Partially supported by the National Science Foundation grant DMS-1266383}
}

\date{\today}

\begin{document}
\maketitle
\pagestyle{plain}
\pagenumbering{arabic}

\begin{abstract}
\no The evolution of a dispersal trait is a classical question in evolutionary ecology, which has been widely studied with several  mathematical models.   The main question is to define the fittest dispersal rate for a population in a bounded domain, and, more recently, for traveling waves in the full space.
\smallskip

\no In the present study, we reformulate the problem in the context  of adaptive evolution. We consider a population structured by space and a genetic trait acting directly on the dispersal (diffusion) rate under the effect of rare mutations on the genetic trait. We show that,  as in simpler models, in the limit of vanishing mutations, the population concentrates on a single trait associated to the lowest dispersal rate. We also explain how to compute the evolution speed towards this evolutionary stable distribution. 
\smallskip

\no The mathematical interest stems from the asymptotic analysis which requires a completely different treatment for each variable. For the space variable, the ellipticity leads to the use the maximum principle and Sobolev-type regularity results.  For  the trait variable, the concentration to a Dirac mass requires  a different treatment.  This is based on the WKB method and viscosity solutions   leading  to an effective Hamiltonian (effective fitness of the population) and a constrained Hamilton-Jacobi equation.
\end{abstract}

\bigskip


%
\section{Evolution of dispersion}
\label{sec:ed}

 \paragraph{Evolutionary dynamics of a structured population.} 
There are several well established mathematical formalisms to describe evolution.  
Game theory is widely used since the  seminal paper \cite{maynard}; see also \cite{HS2}. Dynamical systems are also employed  to describe the possible invasion of a population be a mutant, and to characterize Evolutionary Stable Strategies and other mathematical concepts; see, for instance,~\cite{diekmann, geritz}. Stochastic individual based models are often used to describe the evolution of  individuals undergoing birth, death and mutations. Relations with other approaches can also been made  at  the large population limit, see \cite{BBC15, NC.RF.SM:06}  and the references therein. 
\smallskip

The formalism we use here is still different and relies on a population structured by a phenotypical  trait and competing for a limited resource. This approach was initiated and has been  widely  studied in  \cite{OD.PJ.SM.BP:05}. Several  other versions use the formalism of a population structured by a trait, undergoing mutations and competition \cite{GB.BP:08, DJMR, PJ.GR:09, AL.SM.BP:10}. All these papers, however, consider only proliferative advantage. The population with the highest birth rate or best competition ability survives while other traits undergo extinction. Mathematically this is represented by a limiting process where the population number density, denoted here by $n$,  takes the form of a weighted Dirac concentrated at the fittest trait $\bar \theta$.
\smallskip

The extension to models, where the phenotypical trait is combined  with another structuring variable, usually space, is more recent and leads to considerable mathematical difficulties; see \cite{ADP, EB.SM:14,Turanova_2014, SM_BP, CDM, Campillo_Fritsch}. This is  mainly due to the fact that for the trait variable  $\theta$ the solutions concentrate as described above while in the space variable solutions remain bounded. 
\smallskip

Another motivation for considering the model in this note is to study  the selection of the fittest individuals without a proliferative advantage. In this context, the reproduction rate might be compensated by another advantage. This gives rise to  the question of defining an ``effective fitness.'' The gradient of the effective fitness determines the direction of trait evolution and its maximum defines the evolutionary stable strategy. 
\smallskip

A particularly interesting example, both mathematically and biologically,  in this directions is the selection of a dispersal rate, which we describe next in the context of a continuous dispersal trait.  When only two species are represented by their number densities $n_1(x)$, $n_2(x)$ and are  competing for the same resource (carrying capacity) $K$,  the question is to know which of the two dispersal rates $D_1$, $D_2$ is prefered in the competiton. Then, the model is 
$$
 - D_i \Delta_x n_i  = n_i \big( K(x)- n_1+n_2 \big), \qquad x\in \Omega, \quad i=1, \; 2.
$$
Is it better to `move' faster or slower? In other words, is it favorable to have $D_1 > D_2$ or the contrary?
 
 \paragraph{The model.} Here we assume that all dispersal rates are possible, and include mutations. We ask the question of the rare mutation limit of  the steady state version for a still simple model of evolution of dispersal in a population. The main modeling assumptions are:  (i)~all individuals wear a phenotype characterized by a parameter $\theta$, which induces a dispersal rate $D(\theta)$, (ii)~ a Fisher-type   Lotka-Volterra growth/death rate 
with a space dependent carrying capacity $K$ and limitation by the total population whatever  the trait is, 
and (iii)~rare mutations acting on the genetic variable and modeled by a diffusion with covariance  $\sqrt 2 \e$;
we refer to \cite{NC.RF.SM:06} for a derivation  of this type of equations from individual based stochastic models. 
 
\smallskip

More precisely,  we study the asymptotic behavior, as $\e \to 0$, of the density  $n_\e=n_\e(x,\theta): \Omega \times \R \to (0,\infty)$, with $\Omega\subset \R^d$   a smooth  domain,  of the nonlocal and nonlinear problem 
\beq \begin{cases}
 - D(\theta) \Delta_x n_\e -\e^2 \Delta_\theta   n_\e = n_\e \big( K(x) - \rho_\e(x) \big) \  \  \text{in } \ \ \Omega \times \R, 
 \\[1pt]
\rho_\e(x) = \dis\int_0^1 n_\e(x, \theta) d\theta,
\end{cases}
\label{equation}
\eeq
with Neumann boundary  condition on $\p \Omega $ and $1-$ periodicity in $\theta$, that is,  if $\nu$ is the external normal vector to $\Omega$,
\beq
\f{\p}{\p \nu } n_\e = 0 \ \text{ on} \ \p \Omega\times \R \ \text{and} \ n_\e \ \text{is } \ 1-\text{periodic in} \ \theta.
\label{bc}
\eeq
We have chosen periodic boundary conditions in $\theta$ to simplify some technical aspects concerning a priori estimates. 
\smallskip
 
 As far as the carrying capacity $K$  and the dispersion rate $D$ are concerned, we assume 
\beq
K\in C^1(\bar \Omega),  \qquad   \text{there exists $K_m >0$ such that} \ K(x) \geq K_m >0, \qquad K \text{ is not constant}.
\label{asK}
\eeq
and
\beq \bepa
D \in C^1 \big(\R;(0,\infty) \big)  \  \text{ is} \ 1\text{-periodic and} 
\\[1mm]
\text{ there exists a unique local minimizer  $\theta_m  \in [0,1)$ such that} \   D\geq D(\theta_m  ):=D_m>0; 
\eepa
\label{asD}
\eeq
we note  \eqref{asD} is used to assert that the effective Hamiltonian also has a minimum at $\theta_m$. 
\smallskip

\no  Finally, for technical reasons, we also assume that 
\beq\label{convex}
\Omega \ \text{is convex}.
\eeq

 It follows from, for example, \cite{BNS, CD}, that, given \eqref{asK} and \eqref{asD},  the problem \eqref{equation} and \eqref{bc} admits a strictly positive solution $n_\e:\Omega \times \R \to (0,\infty)$.

\paragraph{Formal derivation of the mathematical result.} We proceed now formally to explain what happens in the limit $\e\to 0$ and, hence, motivate the statement of the results. As it is often the case with problems where it is expected to see concentration in the limit, we make the exponential change of variables 
\beq\label{exp}
n_\e = \exp\big( {u_\e}/{\e} \big),
\eeq
which leads to
\beq
- \f{D(\theta)}{\e}  \Delta_x u_\e -\f {D(\theta)} {\e^2} |\nabla_x u_\e |^2 - \e  \Delta_\theta u_\e - |\nabla_\theta u_\e |^2 = K(x)-\rho_\e (x) \ \text{ in } \ \Omega \times \R,
\label{eqhje1}
\eeq
with 
\beq\label{expbc}
u_\e \ \text{is $1-$periodic in $\theta$} \ \text{ and } \ \f{\p}{\p \nu } u_\e = 0 \ \text{ on} \ \p \Omega\times\R.
\eeq

 It is clear from \eqref{eqhje1} that, if the $u_\e$'s have, as $\e\to 0$,  a limit $u$, it must be independent of $x$, and, hence, it is natural to expect the expansion  
\beq\label{ansatz}
u_\e(x,\theta)=u(\theta) + \e \ln {\mathcal N}(x,\theta) + O(\e).
\eeq

 Assuming that, as $\e\to 0$, the $\rho_\e$'s converge to some $\rho$, a formal computation suggests that $ {\mathcal N}$ is the positive eigenfunction of
\beq\label{ef}
\begin{cases}
- D( \theta ) \Delta_x {\mathcal N}  = {\mathcal N} \;  \big( K - \rho \big) +{ \mathcal N} \; H\big(\theta, \rho(\cdot)\big) \  \text{in} \  \Omega,
\\[1mm]
 \f{\p}{\p \nu }{ \mathcal N} = 0  \ \text{on }  \p \Omega, 
\end{cases}
\eeq
with eigenvalue $H\big(\theta, \rho(\cdot) \big)$, and that $u$ solves the constrained Hamilton-Jacobi equation
\beq\label{chj} 
\begin{cases}
 -|\nabla_\theta u|^2 =  - H\big(\theta,  \rho (\cdot)\big) \ \text{ in } \ \R,
 \\[1mm]
 \max_{\theta \in \R} u(\theta)=0,
 \\[1mm]
u \ \text{is $1-$ periodic}.
\end{cases}
\eeq
The constraint on the $\max u$ becomes evident from the facts that, as it turns out, the $\rho_\e$'s are bounded uniformly in $\e$ and the equality 
\beq
\rho_\e(x)=\int_0^1 n_\e(x,\theta) d\theta =\int_0^1 \exp \big(u_\e(x,\theta)/\e \big) d\theta,
\eeq
which also suggests that, as $\e\to 0$, the $n_\e$'s behave like a Dirac mass with weight $\rho$.

\paragraph{The mathematical result.} To state the result we recall that $\delta$ denotes the Dirac mass at $0$ and we introduce the nonlinear Fisher-type stationary problem
\beq
\begin{cases}
- D_m \Delta_x N_m = N_m \; (K - N_m) \ \text{on} \  \Omega\\[1mm]
 \f{\p}{\p \nu } N_m = 0 \ \text{ on} \ \p \Omega,
\end{cases}
\label{eq:Nm}
\eeq
which, in view of \eqref{asK} and \eqref{asD}, admits a positive solution $N_m>0$ (see, for example, \cite{BNS, CD}.) 

 We have:
\begin{theorem} Assume  \eqref{asK}, \eqref{asD}, and \eqref{convex}. Then,   as $\e \to 0$ and in the sense of distributions, 
$$
n_\e \to N_m(x) \delta(\theta-\theta_m) \quad  \text{ and}  \quad  \rho_\e  \to N_m.
$$
Moreover, as $\e \to 0$ and uniformly in $x$ and $\theta$,  $\e \ln n_\e  \to u$, where $u$ is the unique $1-$periodic solution to \eqref{chj}, with $\rho=N_m$,  such that  $\max u =u(\theta_m)=0.$ 
\label{th1}
\end{theorem}

\paragraph{Biological interpretation.}  
The conclusions of Theorem~\ref{th1} can be thought as a justification of the fact that the population selects the ``slowest''  individuals in accordance with several previous observations on the  evolution of dispersal. In this respect, the eigenvalue $-H\big(\theta, \rho(\cdot) \big)$ defines the fitness of individuals depending upon their trait. This fact can be stated using  the canonical equation \eqref{canonical}, which is formally derived in Section~\ref{sec:conclusion}. In the words of adaptive dynamics,  our result characterizes the unique Evolutionary Stable Distribution (or Strategy), \cite{diekmann,PJ.GR:09}.

\smallskip

That mutants with lower dispersal rates can invade a population, that is the characterization of $\theta_m$ by property~\eqref {asD}, is known from the first mathematical studies \cite{Hast_83,DHMP}. However, these papers use time scale separation, heuristically assuming a mutant appears, with `small mutation', and compete with the resident population. Our approach here is more intrinsic since we consider structured populations competing  for resources and undergoing mutations. Surprisingly when set in the full space where the problem is characterized in terms of traveling waves, the opposite effect is observed, that is mutants with higher dispersal rates are selected giving rise to accelerating waves, \cite{EB.VC.NM.SM:12, EB.SM:14, Turanova_2014, BJS, BMR}. For two competing populations, the combined effect of dispersal and  a drift is studied in \cite{HL_09}.  The analysis of dispersal evolution also gave  rise to the notion of ideal free distribution \cite{CCL_08, CDM}. 

\smallskip

The question of dispersal evolution is a classical and important topic in evolutionary biology. The reader can consult \cite{ronce} for a survey of the many related issues, to \cite{PDGM} for the case with patches and demographic stochasticity, to \cite{HMMV} for the case of trajectories with jumps (nonlocal operators) and to \cite{EB.VC.NM.SM:12} for other biological references about accelerating fronts. A formalism using Fokker-Planck equation is used in \cite{PSL}. 
Also, let us mention that a remarkable qualitative aspect fo spatial sorting  is that in the full space, the largest dispersion rate is selected \cite{EB.VC.NM.SM:12, EB.VC14, EB.VC.GN15}. Finally, another mathematical approach to the concentration effect, stated in Theorem~\ref{th1},  can be found in \cite{LamLou}.

\paragraph{Organization of the paper.} 
In Section~\ref{sec:est_rho} we prove some uniform in $\e$ estimates for the $\rho_\e$'s that  are then used in Section~~\ref{sec:H}  where we derive the effective Hamiltonian, that is the eigenvalue problem of \eqref{ef}.  In Section~~\ref{sec:chj} we introduce  the constrained Hamilton-Jacobi equation to conclude the proof of Theorem~\ref{th1} and Theorem~\ref{th:u}. In Section~~\ref{sec:proofs} we prove the two technical lemmata that were used in Section~~\ref{sec:chj}. Finally,  In Section~~\ref{sec:conclusion}, we provide some perspectives about the problem, namely a more precise asymptotic expansion and the parabolic case, as well as  numerical examples for the evolution driven by the parabolic equation.

\section{Estimates on $\rho_\e$}
\label{sec:est_rho}

We state and prove here some, uniform in $\e$, estimates for the $\rho_\e$'s, which are fundamental for the analysis in the rest of the paper; here $|\Omega|$ is the measure of $\Omega$ and $C(A)$ denotes a constant that depends on $A$.

\begin{lemma} \label{lm:takis2}Assume \eqref{asK} and \eqref{asD}. There exist positive independent of $\e$ constants $C_1=C_1(K,D, \Omega)$, $C_2=C_2(K,D, \Omega)$ and $C_3=C_3(|\Omega|, \min K)$  such that  
\begin{equation}
\label{lm:est_rho}\begin{cases}
(i)~ \ \ 0 \leq  \rho_\e  \leq C_1, \\[1mm]

(ii)~ \  \sup_{\e\in (0,1)} \| \rho_\e\|_{W^{2,p}(\Omega)} \leq C_2   \ \text{ for all $p\in[1,\infty)$, and}\\[1mm]

(iii) ~ \int_\Omega  \rho_\e(x) dx  \geq C_3, \end{cases}
\end{equation}  
and, along subsequences $\e \to 0$,  the $\rho_\e$'s converge uniformly in $\overline \Omega$ to $\rho \in C(\overline \Omega).$ 
\end{lemma}

\proof  We first observe that $\rho_\e$ trivially satisfies the Neumann condition 
\begin{equation}\label{takis1}
\f{\p}{\p \nu } \rho_\e = 0 \ \text{ on} \ \p \Omega.
\end{equation}
After dividing \eqref{equation} by $D(\theta)$, integrating in $\theta$ and using the periodicity in $\theta$,  we find
$$
 - \Delta_x \rho_\e  -\e^2  \int_0^1 n_\e(x,\theta) \Delta_\theta \f 1 D d \theta=  \int_0^1  \f 1 D n_\e d\theta \; (K-  \rho_\e),
$$
and, hence, for some  constant $C$, which only depends on $D$ and $K$, we have
$$
 - \Delta_x \rho_\e +  \f{1}{\| D\|_\infty}   \rho_\e^2 \leq   C  \rho_\e.
$$
Then \eqref{lm:est_rho}(i) follows from the strong maximum principle, while the $W^{2,p}$-estimates are a consequence of  the classical elliptic regularity theory. 

 The lower bound \eqref{lm:est_rho}(iii) comes from integrating  \eqref{equation} in $x$ and $\theta$. Indeed, in view of the assumed  periodicity, we find  
$$
\max K \int_\Omega  \rho_\e(x) dx \geq \int_\Omega  \rho_\e(x) K(x) dx = \int_\Omega  \rho_\e(x)^2 dx \geq \f 1 {|\Omega |}\left(\int_\Omega  \rho_\e(x) dx\right)^2 . 
$$

 The last claim is an immediate consequence of the a priori estimates and the usual Sobolev imbedding theorems.

\hskip6.5in \qed

\section{The effective Hamiltonian}
\label{sec:H}

For  $\theta\in [0,1]$  and $\rho \in L^\infty \big(\Omega; [0,\infty) \big)$ we consider the eigenfunction ${\mathcal N}= {\mathcal N} \big(\cdot; \theta,  \rho(\cdot) \big) \in H^1(\Omega)$ and the eigenvalue 
 $H=H\big(\theta,  \rho(\cdot) \big)$ of 
\beq \begin{cases}
- D( \theta ) \Delta_x {\mathcal N}  = {\mathcal N} \; ( K - \rho ) +{ \mathcal N} \; H \quad  \text{in} \  \Omega,
\\[1mm]
 \f{\p}{\p \nu }{ \mathcal N} = 0  \ \text{on }  \p \Omega  \qquad \text{and} \ \   \int_\Omega {\mathcal N}  \big(x; \theta, \rho(\cdot) \big)^2 dx =1.
 \end{cases}
\label{effective_hamiltionian}
\eeq
Note that in view of  \eqref{asK} and the regularity of $\Omega$, the existence of the pair $ \big({\mathcal N}  \big(\cdot; \theta,  \rho(\cdot) \big), H \big(\theta,  \rho(\cdot)  \big)  \big)$ follows from, for example, \cite{BNS, CD}.
\smallskip

 The next lemma provides some important estimates and information about $H \big(\theta,  \rho(\cdot)  \big)$. In the statement we use the notation $K_M:=\max_{\overline \Omega} K.$
\begin{lemma} Assume \eqref{asK} and \eqref{asD}. Then
\beq
(i) \  -K_M \leq H\big(\theta,  \rho(\cdot)\big) \ \  \text{and} \ \  (ii) \  H \big(\theta,  \rho(\cdot)\big)  \leq \f{1}{|\Omega |} \int_\Omega \rho(x) dx, 
\label{H:est1}
\eeq
\beq
\text{ the maps $\theta \to H\big(\theta , \rho(\cdot)\big)$ and $\theta\to  D(\theta)$ have the same monotonicity properties for all $\rho$}, 
\label{H:est3}
\eeq
and, in particular, $\theta_m$ is the unique minimum of $H(\theta, \rho)$  in $[0,1]$  for all $\rho$. 
\label{lm:H}
\end{lemma}
\proof Multiplying the equation by  ${\mathcal N}$ and integrating in $x$ gives 
$$
0 \leq D(\theta) \int_\Omega |\nabla {\mathcal N} |^2  = \int_\Omega {\mathcal N}^2 [K-\rho +H] dx \leq  \int_\Omega {\mathcal N}^2 dx \;  [K_M + H], 
$$
and, thus,   \eqref{H:est1}(i)  holds. 
\smallskip

  The upper bound  \eqref{H:est1}(ii) follows from the positivity of $K$, since,  after dividing the equation by $ {\mathcal N}$ and integrating  by parts, we find
$$
- \int_\Omega \f{\big| \nabla_x  {\mathcal N} \big|^2}{  {\mathcal N}^2  } dx =   \int_\Omega   [ K -\rho + H] \, dx  \leq 0. 
$$

\smallskip

 For  \eqref{H:est3}, we differentiate in $\theta$ the equation in \eqref{effective_hamiltionian} to find 
\beq\label{takis}
- D'(\theta) \Delta_x {\mathcal N} - D(\theta) \Delta_x {\mathcal N}_\theta =  {\mathcal N}_\theta \;  \big( K(x) - \rho(x) + H\big) +  {\mathcal N} \;  H_\theta,
\eeq
where  ${\mathcal N}_\theta$  and  $H_\theta$ denote derivatives with respect to $\theta$, we multiply by $ {\mathcal N} $ and  integrate by parts using the 
boundary condition to get
\beq\label{takis2}
D'(\theta) \int_\Omega \big| \nabla_x  {\mathcal N} \big|^2 + D(\theta) \int \nabla_x  {\mathcal N}.\nabla_x {\mathcal N}_\theta   =
\int  {\mathcal N}  {\mathcal N}_\theta [K-\rho]+   H_\theta \big(\theta; \rho(\cdot) \big). 
\eeq
Next we use the fact that  the $L^2-$normalization of $ {\mathcal N}$ yields $\int_\Omega  {\mathcal N}_\theta {\mathcal N} dx =0$ and after multiplying  the equation in \eqref{effective_hamiltionian} by  $ {\mathcal N}_\theta$,  integrating by parts and subtracting  the result from \eqref{takis2}  we find 
$$
D'(\theta) \int_\Omega \big| \nabla_x  {\mathcal N} \big|^2 = H_\theta \big(\theta; \rho(\cdot) \big). 
$$ 
Since  we have assumed in \eqref{asK} that $K$ is not constant, $\nabla_x  {\mathcal N}$ does not vanish and the result follows.
\smallskip

 To conclude, we observe that  \eqref{asD} and \eqref{H:est3} yield that $ H$  as the same monotonicity, in $\theta$,  as $D$ and, thus,  has a unique local minimum at $\theta_m$. 

\hskip6.3in \qed

\section{The constrained Hamilton-Jacobi equation}
\label{sec:chj}

We prove here the generalized Gaussian-type convergence asserted in Theorem~\ref{th1}, derive the constrained Hamilton-Jacobi equation \eqref{chj} and state some more properties. For the benefit of the reader we restate these assertions as a separate theorem below.

\begin{theorem}  The family  $u_\e$  is uniformly in $\e$ Lipschitz continuous and converges, uniformly in $x$ and $\theta$, to $u$,  which is independent of $x$ and satisfies, in the viscosity sense, the constrained Hamilton-Jacobi equation \eqref{chj}.  Moreover, 
\beq
\min_\theta H\big(\theta, \rho(\cdot) \big) =0= H\big(\theta_m, \rho(\cdot) \big)=  H_\theta \big(\theta_m, \rho(\cdot) \big).
\label{eq:ESS}
\eeq
\label{th:u}
\end{theorem}
Since $-H$ represents the fitness, it turns out that  \eqref{eq:ESS} characterizes $\delta (\theta - \theta_m)$ as the Evolutionary Stable Distribution (or $\theta_m$ as the Evolutionary Stable Strategy). See \cite{diekmann,PJ.GR:09}.

It follows  from \eqref{eq:ESS} that both $H$ and its derivative vanish at $\theta_m$. As a result  the viscosity solution $u$ of~\eqref{chj} also vanishes at $\theta_m$ as $u(\theta) =O\big((\theta- \theta_m)^{3/2}\big)$, and this makes the connection with the result of \cite{LamLou}.

\smallskip

 The proof is a consequence of the next two  lemmata which are proved later in the paper.

\begin{lemma} [Bounds on $u_\e$] There exists an independent of $\e$ constant $C$ such that 
\beq
\int_\Omega \max_\theta u_\e (x, \theta) dx \leq C \e .
\label{eq:u_upper}
\eeq
\label{lm:ubdd}
\end{lemma}
\begin{lemma} [Lipschitz estimates] There exist  an independent of $\e$ constant $C$ such that 
$$
\f1{\e^2}  |\nabla_x u_\e |^2 +  |\nabla_\theta u_\e |^2 \leq C \quad \text{and} \quad  \max_{x \in \Omega, 0 \leq \theta\leq 1}  u_\e (x, \theta) \leq C\e.
$$
Moreover, the $u_\e$'s converge, along a subsequence $\e \to 0$ and uniformly in $x$ and $\theta$, to a Lipschitz  and $1$-periodic function $u:\R\to \R$ such that
$
\max_{ 0\leq \theta \leq 1 } u ( \theta)  = 0.
$
\label{lm:ulip}
\end{lemma}

 We continue with the proofs of Theorem~\ref{th:u} and Theorem~\ref{th1}.
\smallskip

\noindent{\bf Proof of Theorem \ref{th:u}.} 
The fact  that any limit $u$ of the $u_\e$'s satisfies \eqref{chj}  is a standard application of the so-called perturbed test function method (see \cite{evans}) and we do not repeat the argument. 
\smallskip

 It follows from \eqref{chj} that  $H\big(\theta, \rho(\cdot) \big) \geq 0$,  while, at any maximum point   $\bar {\theta} $ of $u$, we must have $H\big(\bar \theta, \rho(\cdot) \big) \leq 0$, and, hence,  $H\big(\bar \theta, \rho(\cdot) \big) = 0,$ and, in view of Lemma~\ref{lm:H}, 
$$
\min_\theta H\big(\theta, \rho(\cdot) \big) =0= H\big(\theta_m, \rho(\cdot) \big).
$$
As a result the only possible maximum point of any solution of \eqref{chj} must be at at $\theta_m$, which implies that the equation has a unique solution. 

\smallskip

 Also the knowledge of $\theta_m$ determines uniquely the limit  $\rho(x) =  N_m(x)$, from equation \eqref{eq:Nm}, and, thus,  the full family $(\rho_\e, u_\e)$ converges. 

\hskip6.3in \qed

\smallskip

\noindent{\bf Proof of Theorem \ref{th1}.}  The statement in terms of $n_\e$ is an immediate consequence of Theorem~\ref{th:u}. Because $u_\e$ achieves a unique maximum at $\theta_m$, from the Laplace formula for $n_\e$ written as \eqref{exp}, we conclude that the $n_\e(x, \theta)$'s  converge weakly in the sense of measures to $\rho(x) \delta(\theta - \theta_m)$, with $\rho(x)$ the limit of $\rho_\e$ (see Section~\ref{sec:est_rho}).

\no Next, integrating equation \eqref{equation} in $\theta$ we conclude that
$$
- \Delta \int_0^1 D(\theta) n_\e d\theta = \rho_\e(x) \big( K(x) - \rho_\e(x) \big). 
$$
Passing to the limit $\e \to 0$,  and taking into account that $n$ is a Dirac mass at $\theta_m$, we find that $\rho =N_m$ because they both satisfy the equation \eqref{eq:Nm}.

\hskip6.3in \qed

\section{The proofs}
\label{sec:proofs}

We begin with the proof of  Lemma~\ref{lm:ubdd}
\smallskip

\noindent{\bf Proof of Lemma \ref{lm:ubdd}.} Integrating \eqref{eqhje1}, we find 
\beq
\f 1 {\e^2}\int \int  |\nabla_x u_\e |^2 dx d\theta  + \int \int  |\nabla_\theta u_\e |^2 dx d\theta \leq C. 
\label{hj:h1est}
\eeq
 The claim for the maximum bound follows from  the $L^\infty-$estimate on $\rho_\e$.  
\smallskip 
 
 For each $x\in \Omega$, let $M_\e (x): = \max_\theta  u_\e(x, \theta)$ and choose $\theta_\e(x)$ such that  $M_\e (x)=u_\e (x, \theta_\e(x)).$ Then  it follows from  \eqref{lm:est_rho}(i) that 
\beq
 e^{M_\e(x)/\e}   \int_0^1 e^{\big(u_\e(x, \theta)-u_\e(x, \theta_\e)\big)/\e} d\theta =\rho_\e (x) \leq C.
\label{hjc:estu}
\eeq
Inserting in \eqref{hjc:estu} the  estimate
$$
u_\e(x, \theta) - u_\e(x, \theta_\e(x)) = \int_\theta^{\theta_\e(x)} \p_\theta u_\e(x, \theta') d\theta' \geq - 1 - \int_0^1 | \p_\theta u_\e(x, \theta')|^2 d\theta',
$$
we get  
$$
e^{M_\e(x)/\e} \; e^{\dis - \f 1 \e [1+\int_0^1 | \p_\theta u_\e(x, \theta')|^2 d \theta'] }  \leq C,
$$
and, hence,
$$
\f {M_\e(x)}{\e} -\f 1 \e[1+  \int_0^1 | \p_\theta u_\e(x, \theta')|^2 d \theta'] \leq \ln C.
$$
Using \eqref{hj:h1est}, we obtain 
$$
\int_\Omega M_\e(x) dx \leq C \e,
$$
and this concludes the proof.

\hskip6.3in \qed

\smallskip
Next we discuss the proof of Lemma~\ref{lm:ulip}.

\noindent{\bf Proof of Lemma \ref{lm:ulip}.} We first assume  the Lipschitz  bound and prove the rest of the claims.

\smallskip

 Let ${\mathcal M_\e}:= \max_{x\in\overline \Omega} M_\e (x)$, with $M_\e (x)$  as in the proof of the previous lemma. It is immediate from \eqref{eq:u_upper} and the fact that $|\nabla_x u_\e |^2 \leq \e C$ that, for some $C>0$,
$$
{\mathcal M_\e} \leq \e C.
$$

 Next we show that  $\liminf_{\e \to 0} {\mathcal M_\e} \geq 0$, which, in turn, yields that $\lim_{\e\to 0}{\mathcal M_\e} =0$.
\smallskip

 Let $y_\e \in \overline \Omega$ be such that 
${\mathcal M_\e}:= M_\e(y_\e) = u_\e \big(y_\e, \theta_\e(y_\e)\big)$  and write   
$$
\rho_\e (x) = e^{\mathcal M_\e/\e}   \int_0^1 e^{\big(u_\e(x, \theta)-u_\e(y_\e, \theta_\e)\big)/\e} d\theta.
$$
 Combining the lower bound  on $\rho_\e$ in  Lemma~\ref{lm:takis2} and the (Lipschitz) estimate $|\nabla_x u_\e |^2 \leq \e C$,
we get    
$$
C_3 \leq  \int_\Omega \rho_\e (x) )\leq  |\Omega|  e^{\mathcal M_\e/\e}   e^{C}  \int_0^1 e^{\big(u_\e(y_\e, \theta)-u_\e(y_\e, \theta_\e)\big)/\e} d\theta \leq   e^{\mathcal M_\e/\e}   e^{C}  |\Omega| ,
$$
and, thus,  $ \mathcal M_\e \geq - C \e$.
\smallskip

 Now we turn to the proof of the Lipschitz bounds, which is an appropriate modification of the classical  Bernstein estimates to take into account the different scales. We note 
 and prove Lemma~\ref{lm:ulip}. Note 
 that the convexity assumption on $\Omega$ is  used solely  in this proof. 
\smallskip

 We begin by writing the equations satisfied by $|D_x u_\e|^2$ and $ |\p_\theta u_\e|^2$ which we obtain by differentiating \eqref{eqhje1} in $x$  and $\theta$ 
and multiplying by $D_x u_\e$ and  $\p_\theta u_\e$. We have: 
\beq\bea
- \f{D(\theta)}{\e} \Delta_x |D_x u_\e|^2 - &\e \Delta_\theta |D_x u_\e|^2+ 2 \f{D(\theta)}{\e} |D^2_{xx}Ê u_\e|^2 +2 \e  |D^2_{x\theta}Ê u_\e|^2
\\[15pt]
&-2 \f{D(\theta)}{\e^2} D_x u_\e. D_x |D_x u_\e|^2 
- 2 \p_\theta u_\e. \p_\theta  |D_x u_\e|^2= 2 D_x (K- \rho_\e) . Du_\e, 
\eea 
\label{ge:ux}
\eeq
and 
\beq\bea
- \f{D(\theta)}{\e}&\Delta_x |\p_\theta u_\e|^2 - \e \Delta_\theta |\p_\theta u_\e|^2+ 2 \f{D(\theta)}{\e} |D^2_{x\theta}Ê u_\e|^2 +2 \e  |D^2_{\theta \theta}Ê u_\e|^2
\\[15pt]
& - 2 \f{D(\theta)}{\e^2} D_x u_\e. D_x |\p_\theta u_\e|^2 
- 2 \p_\theta u_\e. \p_\theta  |\p_\theta u_\e|^2 = 2  \f{D'(\theta)}{\e} \Delta_x u_\e \p_\theta u_\e+ 4  \f{D'(\theta)}{\e^2} |D_x u_\e|^2 \p_\theta u_\e . 
\eea 
\label{ge:utheta}
\eeq
Let  $q= \f{|D_x u_\e|^2}{\e^2} +  |\p_\theta u_\e|^2$ and compute
$$
\bea
-& \f{D(\theta)}{\e} \Delta_x q - \e \Delta_\theta q+ 2 \f{D(\theta)}{\e^3} |D^2_{xx}Ê u_\e|^2 + \f{1+D(\theta)}{\e}  |D^2_{x\theta}Ê u_\e|^2
+2 \e  |D^2_{\theta \theta}Ê u_\e|^2
\\[15pt]
&-2 \f{D(\theta)}{\e^2} D_x u_\e. D_x q
- 2 \p_\theta u_\e. \p_\theta q = \f{2}{\e^2}  D_x (K- \rho_\e) . Du_\e + 2  \f{D'(\theta)}{\e} \Delta_x u_\e  \p_\theta u_\e+ 4  \f{D'(\theta)}{\e^2} |D_x u_\e|^2 \p_\theta u_\e . 
\eea 
$$
\\
Assume that   $(\bar x, \bar \theta)$ is a maximum point of $q$. Because of the convexity assumption, we have $\f{\p}{\p \nu}Êq <0$ on the boundary (see \cite{PLL_1985}) and thus  $x_\e \notin \p\Omega \times [0,1]$.
\smallskip

 Therefore, at this point $(\bar x, \bar \theta)$, we have
$$\bea
 \f{D}{\e^3} |D^2_{xx}Ê u_\e|^2 + \e  |D^2_{\theta \theta}Ê u_\e|^2
 \leq \f{1}{\e^2}  D_x (K- \rho_\e) . Du_\e +   \f{D'}{\e} \Delta_x u_\e  \p_\theta u_\e+ 2  \f{D'}{\e^2} |D_x u_\e|^2 \p_\theta u_\e 
 \\[15pt]
  \leq \f{1}{\e^2}  D_x (K- \rho_\e) . Du_\e +  \delta  \f{D}{\e^3}( \Delta_x u_\e )^2 +\f{\e}{\delta} (  \p_\theta u_\e)^2+  2  \f{D'}{\e^2} |D_x u_\e|^2 \p_\theta u_\e 
 \eea 
$$
 and we choose $\delta$ small enough so that we can absorb the term $ \delta  \f{D}{\e^3}( \Delta_x u_\e )^2$ in the left hand side.

 Since there is a constant $C_4(d, D)$ such that
$$
\f 1 \e \left[Ê  \f{D}{\e} \Delta_x u_\e + \e \Delta_\theta u_\e \right]^2  \leq C_1 \f{D}{\e^3} |D^2_{xx}Ê u_\e|^2 + C_1 \e  |D^2_{\theta \theta}Ê u_\e|^2,
$$
we conclude (using  the equation) that, for some $C_5(d, D, \delta)>0$, 
$$
\f 1 \e \left[  D \f{|D_x u_\e|^2}{\e^2} +  |\p_\theta u_\e|^2 +K-\rho_\e \right]^2  \leq C_2\left[ \f{1}{\e^2}  D_x (K- \rho_\e) . Du_\e +  {\e} ( \p_\theta u_\e)^2+   \f{2}{\e^2} |D_x u_\e|^2  \p_\theta u_\e \right]. 
$$
 It follows that there exists some positive constant $C$ such that 
$$\bea
q^2 &\leq C \left[1 + \f{1}{\e} | Du_\e | +\e^2  (\p_\theta u_\e)^2+   \f{2}{\e} |D_x u_\e|^2 |\p_\theta u_\e |\right]
\\[15pt]
&\leq C \left[ 1+ \sqrt q + \e^2 q +  \e q^{3/2}\right] .
\eea
$$
From this  we conclude that $q$ is bounded and the Lipschitz continuity statement  is proved. 

\hskip6.3in \qed

%
\section{Conclusion and perspectives}
\label{sec:conclusion}

\paragraph{Conclusion} We have studied a steady state model describing the Evolutionary Stable Distribution for a simple model of dispersal evolution. The novelty is that mutations acting on the continuous dispersal trait (the diffusion rate) are modeled thanks to a Laplacian operator, and this replaces the standard `invasion of a favorable mutant' in the usual time scale separation approach. When the mutation rate is small, we have shown that the minimal dispersal is achieved, in accordance with previous analyses. This is mathematically stated as the limit to a Dirac mass which selects the minimum of the diffusion coefficient in the equation. Technical difficulties rely on a priori estimates in order to make the approach rigorous and establish the constrained Hamilton-Jacobi equation which defines the potential in the Gaussian-like concentration.
\\

We indicate two possible extensions of our results. The first concerns the way to make more precise the convergence result of Theorem~\ref{th1}. The second is about the time evolution problem.
\paragraph{A more precise convergence result.} The question we address here is whether it is possible to make more precise the convergence $n_\e \to N_m(x) \delta(\theta-\theta_m) $ in the weak sense of measures stated  in Theorem~\ref{th1}. 

\smallskip 

The gradient bound in Lemma~\ref{lm:ulip} implies that 
$\| u_\e(x, \theta) -u_\e(0, \theta)\|_{\infty}Ê\leq C \| \nabla u_\e(x, \theta \|_{\infty} \leq C \e$, and,  therefore along subsequences,  the  $\f{u_\e(x, \theta) -u_\e(0, \theta)}{\e} $'s converge  in  $L^\infty\! \!-\  \!w \star$. To prove more about the corrector, it necessary to have an  estimate for $u_\e(0, \theta) - u(\theta)$,  which is not, however,  available from what we have here.
\smallskip

Another approach is to introduce, instead of $u(\theta)$, the  eigenvalue problem 
$$
-\e^2 \Delta_\theta W_\e = W_\e [ - H(\theta, \rho_\e) + \lb_\e ], \quad W_\e \text{ is $1$-periodic}, \qquad W_\e >0.
$$
The change of unknown $w_\e = \e \ln W_\e$ gives the equation 
$$
-\e \Delta_\theta w_\e -|D_\theta w_\e|^2  = - H(\theta, \rho_\e) + \lb_\e .
$$
Standard gradient estimates yield that  $|D_\theta w_\e|$ is bounded independently of $\e$, and, from the equation, we see that the $\e \Delta_\theta w_\e$'s are also bounded independently of $\e$. 
It follows that, as $\e \to 0$,  $\lb_\e \to 0$ and, after an appropriate normalization, $w_\e \to u$ uniformly.

As a result we can factorize the solution of \eqref{equation} as 
$$
n_\e(x, \theta) = v_\e (x, \theta) e^{w_\e(\theta)/\e},
$$
and we claim that, for some other factor $\bar \rho$,
\beq \label{precise}
v_\e(x,\theta) \to \bar \rho \mathcal{N}(x, \theta), \quad  L^\infty\! \!- \!w \ \star, 
\eeq
a statement which is more precise than that of Theorem~\ref{th1}. 

To see this, we write the equation for $v_\e$ as 
$$
-\Delta_x v_\e  - \e^2 \Delta_\theta v_\e -2 \e D_\theta w_\e \cdot D_\theta v_\e = v_\e [K-\rho_\e  + H(\theta, \rho_\e) - \lb_\e ].
$$
 
To obtain bounds on $v_\e$ we notice that
$ \rho_\e (x) = \int_0^1 v_\e W_\e d \theta $ is bounded (from above and below), and,  since $ W_\e \to \delta (\theta -\theta_m)$, we conclude that, for some  $\bar \theta_\e$ and $\underline \theta_\e$ near $\theta_m$ and independent of $\e$ constants $C$ and $c$,  $v_\e (x, \bar \theta_\e) \leq C$ and $v_\e (x, \underline \theta_\e) \geq c$. It then follows from standard arguments that the $v_\e$'s  are  bounded in $L^\infty$.
\smallskip

The limiting equation is  the eigenfunction problem~\eqref{effective_hamiltionian} and positive solutions are all proportional to $\mathcal{N}(x, \theta)$, which gives the statement.

\paragraph{The parabolic problem.} Our approach is mainly motivated by the dynamics of the evolution of dispersal. The steady states are the Evolutionary Stable Distribution and are obtained as the long time distribution of competing populations \cite{PJ.GR:09}. This leads to the study of the time dependent problem
\beq \begin{cases}
\e n_{\e,t}(x, \theta, t) - D( \theta) \Delta_x n_\e -\e^2 \Delta_\theta= n_\e \big( K(x) - \rho_\e(x,t) \big) \ \text{in} \  \Omega\times \R \times (0,\infty), 
\\[1mm]
\rho_\e(x,t) = \int_0^1 n_\e(x, \theta, t) d\theta,
\end{cases}
\label{eq:parabolic}
\eeq
with the Neuman boundary conditions on $\p \Omega$ and $1-$periodicity in $\theta$. 
\smallskip

For this problem there are  two limits of interest, namely $\e \to 0$ and $t \to \infty$. So far, we have studied the limit $t\to \infty$, reaching a steady state~\eqref{equation} of \eqref{eq:parabolic}, and then considered  the limit $\e \to 0$. 
\smallskip

Reversing the order, we need to study first what happens as $\e \to 0$. In this case, we expect that, in the weak sense of measures, 
$$
n_\e(x, \theta, t) \to  {\bar N}(x,t) \delta(\theta - \bar \theta (t)) \quad \text{and}Ê\quad  \rho_\e \to {\bar N}(x,t),
$$
where at least formally, the weight ${\bar N}(x,t) $  is defined, for each $t$,   by the stationary Fisher/KPP equation
\beq\begin{cases}
 - D(\bar \theta(t)) \Delta_x  {\bar N}  =  {\bar N} \big( K(x) -  {\bar N } \big) \ \text{in}  \ \Omega, 
 \\[1mm]
 \f{\p}{\p \nu }  {\bar N}(x,t)  = 0 \ \text{ on} \  \p \Omega.
\end{cases}
\label{eq:No}\eeq 
We can follow the same derivation as before, and discover that the value $\bar \theta(t)$ of the fittest  dispersal trait is now obtained through the time evolution constrained Hamilton-Jacobi equation 
\beq \begin{cases}
u_t(\theta, t) - D(\theta)  |\nabla_\theta u |^2 = - H\big(\theta,  {\bar N}(\cdot, t) \big) \ \text{in } \R,
\\[1mm] 
 \max_\theta u(\theta, t) =0= u(\bar \theta(t),t ) .
\end{cases}
\label{eq:u}
\eeq
Here the effective fitness $H\big(\theta,  {\bar N}(\cdot, t) \big)$ is still given by the eigenvalue problem \eqref{effective_hamiltionian} with $\rho = \bar N$.

\smallskip

Note that, since derivatives  vanish at a maximum point,  we conclude that
$$
H\big( \bar \theta (t),  {\bar N}(\cdot, t) \big) =0,
$$
and we also have $\bar N(x,t)= {\mathcal N}\big(x, \bar \theta(t), \bar N(\cdot,t)\big)$, where   $\bar N(x,t)= \lim_{\e \to 0} \rho_\e(x,t)$,

\smallskip

We recall that, still formally, we can derive from \eqref{eq:u} a  canonical equation for the fittest trait $\bar \theta(t)$ which takes the form, see \cite{OD.PJ.SM.BP:05, AL.SM.BP:10, SM_JMR_cr15, SM_JMR_15}, 
\beq
\f{d}{dt} \bar \theta(t) = - \big( -D_\theta^2 u(\bar \theta(t), t) \big)^{-1}. \nabla_\theta  H \big(\bar \theta (t),  {\mathcal N}(\cdot, t) \big).
\label{canonical}
\eeq
We note that this equation also describes the fact that $\bar \theta (t)$ will evolve towards smaller values of $H$ and thus of the dispersal rate $D$.

The main difficulties compared to the stationary case are to derive a priori estimates for $\rho_\e (x,t)$ analogous to those in Lemma~\ref{lm:takis2} and to obtain gradient estimates on $u_\e$. Since the Lipschitz regularity of  $u$ is optimal and only differentiability can be proved at the maximum point  \cite{GB.BP:08}, giving a meaning to 
 \eqref{canonical} is also a challenge.
 
\begin{figure}[h]
	\centering
		\includegraphics[width=.35\textwidth]{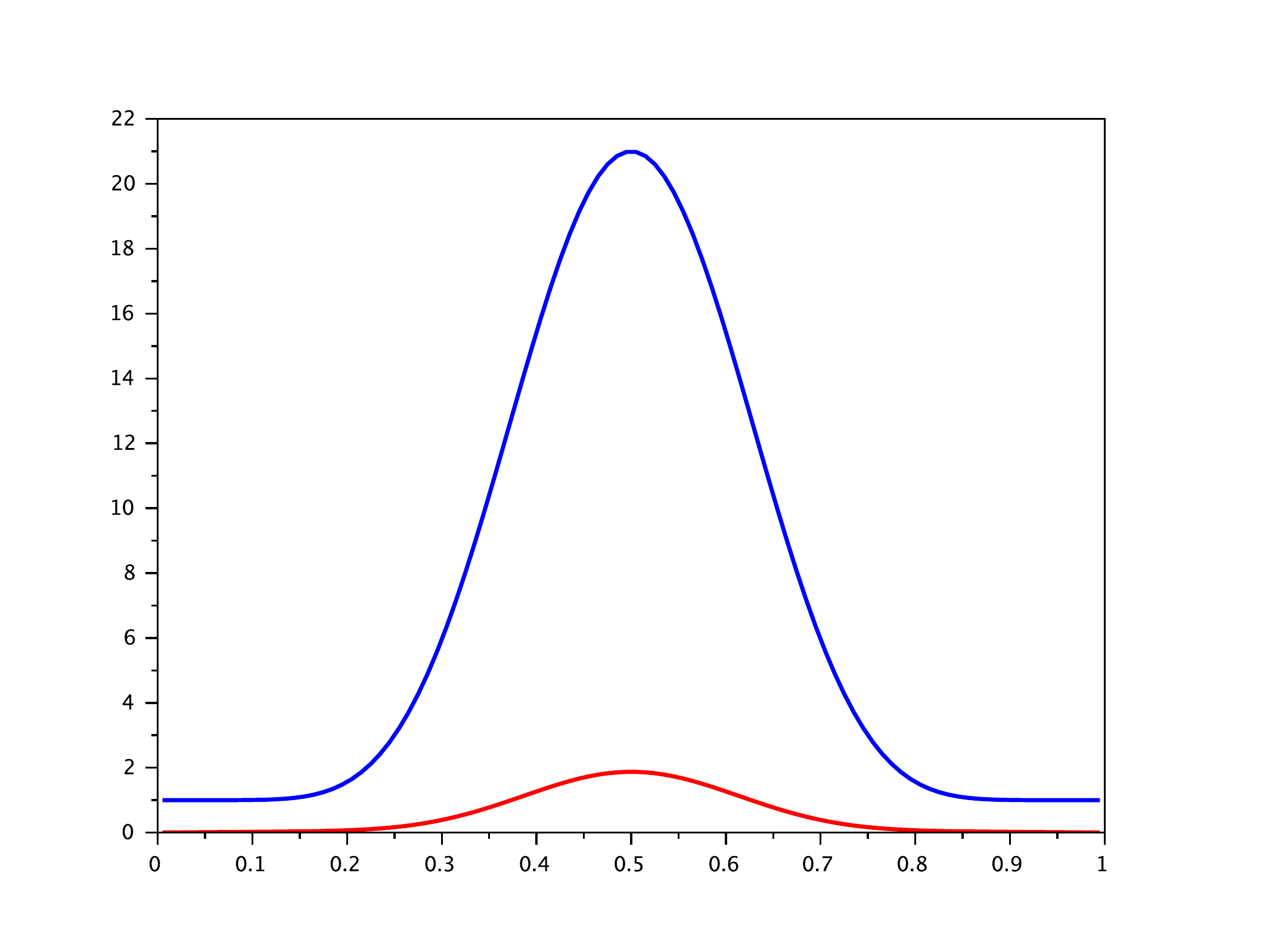}\includegraphics[width=.35\textwidth]{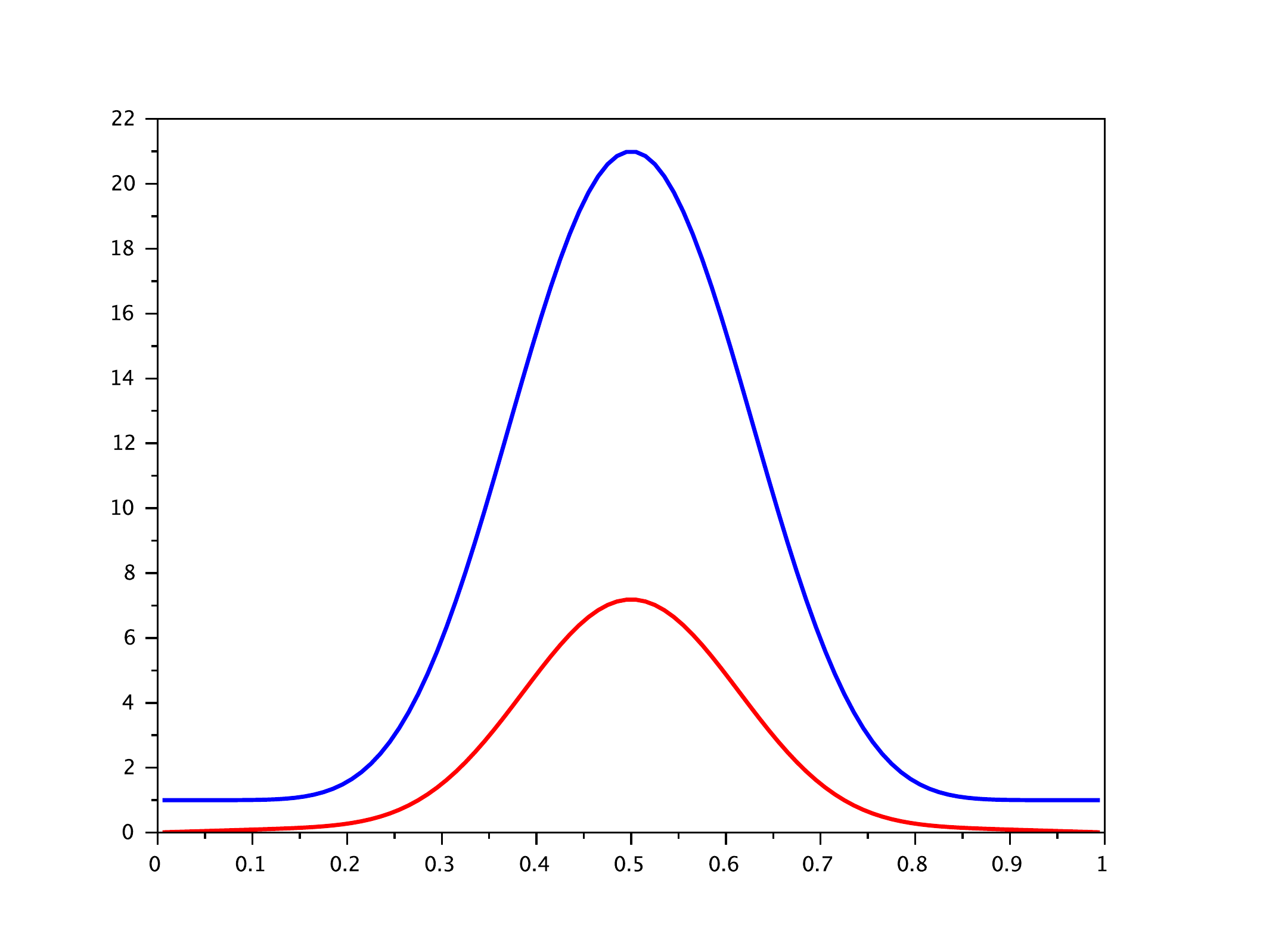}\includegraphics[width=.35\textwidth]{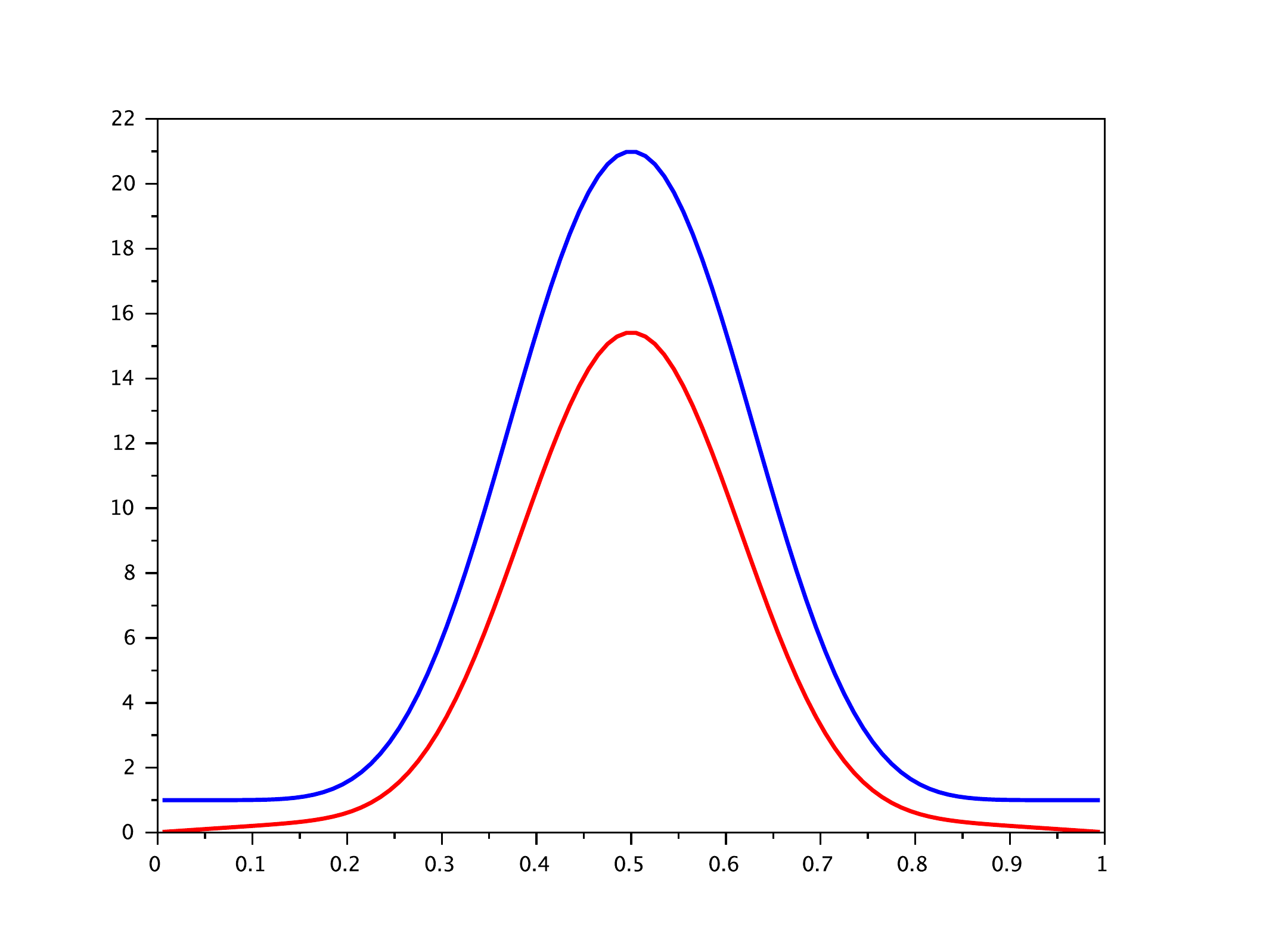}
		
\vspace{-5mm}	
	
		\includegraphics[width=.35\textwidth]{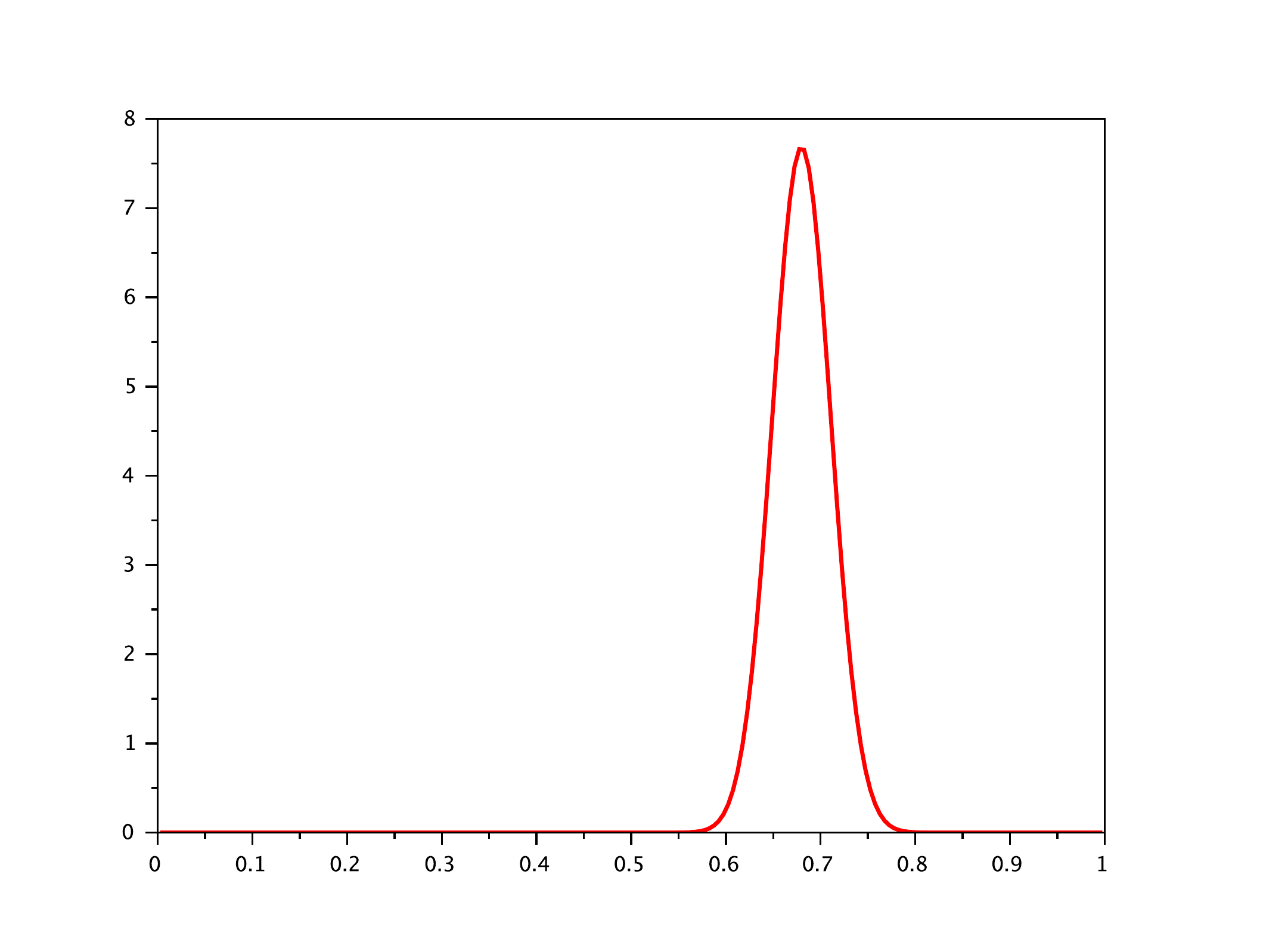}\includegraphics[width=.35\textwidth]{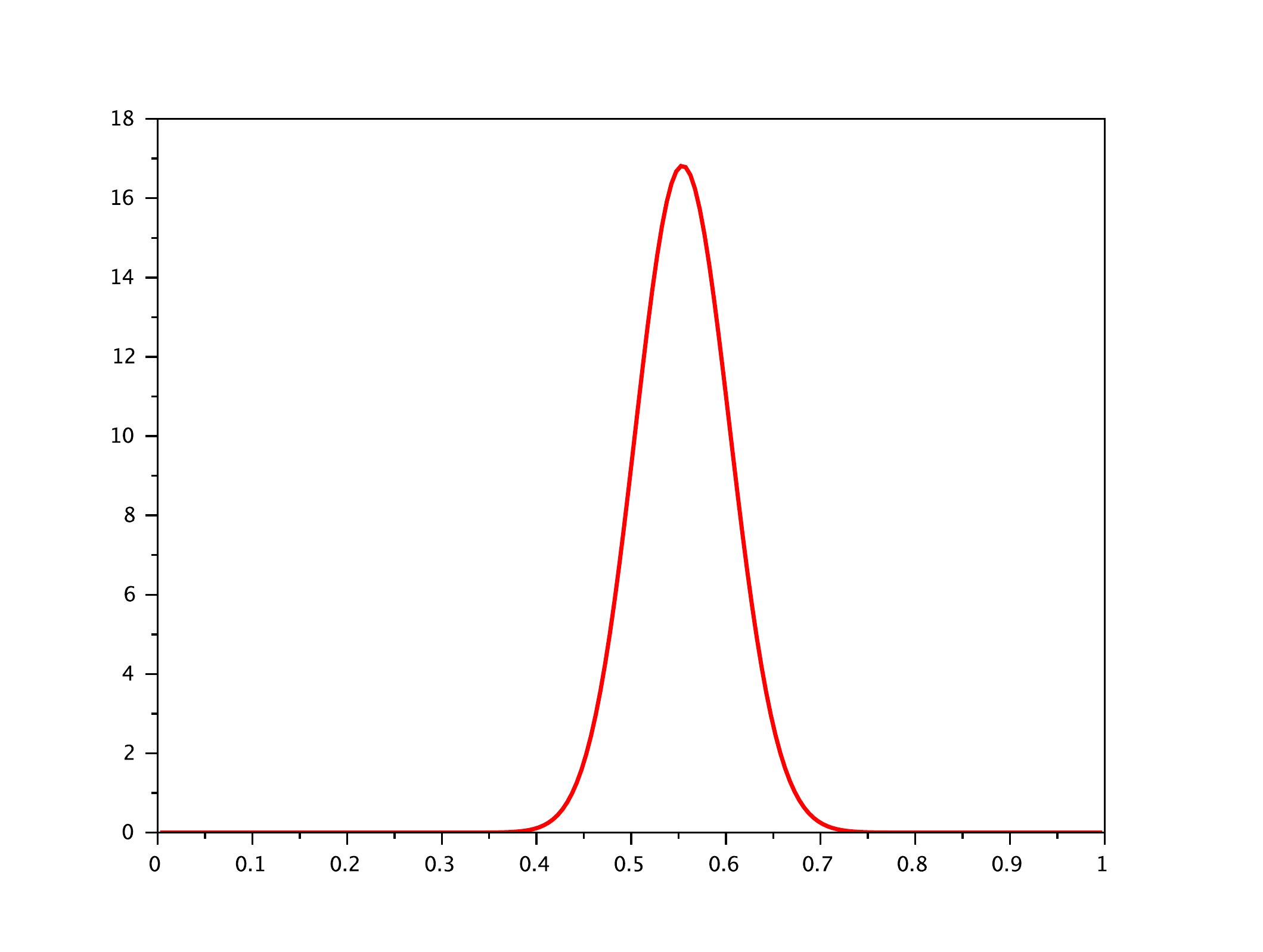}\includegraphics[width=.35\textwidth]{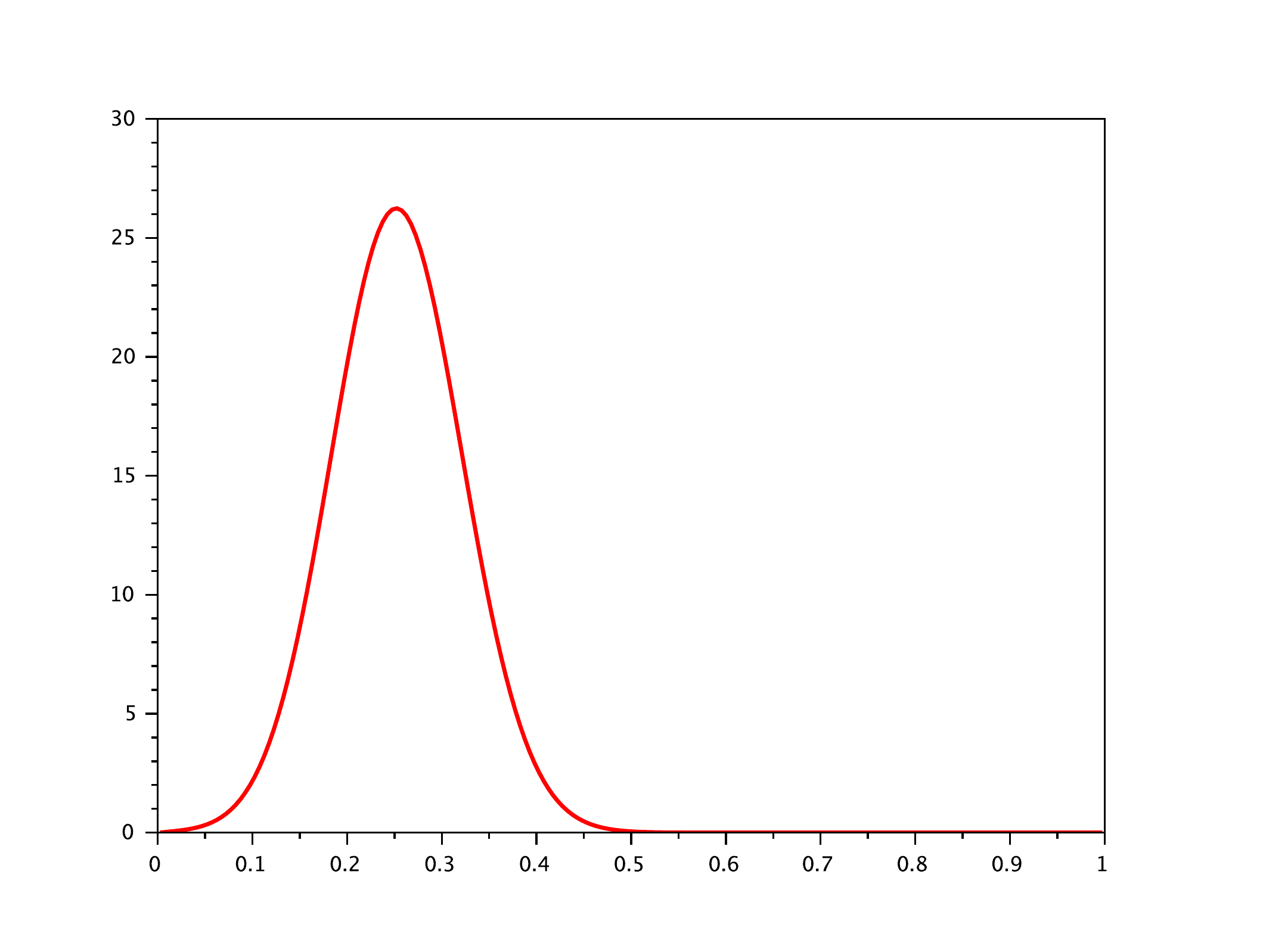}
\vspace{-10mm}		
	\caption{Snapshots of the evolution of the dispersal for three different times according to the parabolic problem~\eqref{eq:parabolic}. Top figures depict the space repartition $\rho(x)$ (fixed blue curve = $K$, increasing red = $\rho$). Bottom figures depict the trait distribution $\int_0^1 n(x, \theta, t) dx$. }
	\label{fig1}
\end{figure}

\paragraph{Some numerics on the parabolic problem.} 
 Numerical simulations are presented in Figure~\ref{fig1} which illustrates the selection of lowest dispersal rate. Considering the parabolic problem~\eqref{eq:parabolic},  that means the convergence of the fittest trait $\bar \theta(t)$ to the smaller values of  the trait (with the coefficients below, this value is $\theta=0$) as $t \to \infty$. 
 
For this simulation we have chosen  $\Omega= (0,1)$ and the data
$$
D(\theta)= 1.5 \; \theta, \qquad K(x) = 1+20 \big( 1-4(x-.5)^2 \big)^8,
$$ 
and we have used,  for the convenience of numerics, Dirichlet boundary conditions both in $x$ and $\theta$. The numerical scheme is the standard three points scheme in each direction, implicit in $x$ and explicit in $\theta$ because $\e= 10^{-2}$ is small enough so as not to penalize the computational time.

In Figure~\ref{fig1}, we observe that the  average trait, which initially is $\theta \sim 0.7$, decreases and gets close to $0$ at the third time displayed here, also the trait  exhibits a concentrated pattern around its average. 

%
%
%

\end{document}